\documentclass[11pt,reqno]{amsart}
\usepackage{amsmath,amssymb,amscd,amsfonts}
\usepackage{enumerate}
\newtheorem{thm}{Theorem}[section]
\newtheorem{lem}[thm]{Lemma}

\newtheorem{prop}[thm]{Proposition}
\newtheorem{cor}[thm]{Corollary}

\newtheorem{assu}[thm]{Assumption}

\newcommand{\Pic}{\operatorname{Pic}}

\newcommand{\OO}{{\mathcal {O}}}

\newcommand{\OX}{{\mathcal{O} _X}}

\newcommand{\Si}{\Sigma}

\newcommand{\Q}{{\mathbb Q}}
\newcommand{\N}{{\mathbb N}}

\newcommand{\F}{{\mathbb F}}
\newcommand{\pp}{{\mathbb{P}}}

\title[Even sets of four nodes on rational surfaces]{%
Even sets of four nodes  on rational surfaces}

\author[A.\ Calabri, C.\ Ciliberto and M.\ Mendes Lopes]{%
Alberto Calabri, Ciro Ciliberto and Margarida Mendes Lopes}



\thanks{The present collaboration takes place in the framework
of the european contract  EAGER, no. HPRN-CT-2000-00099. The first
two authors are members of G.N.S.A.G.A.-I.N.d.A.M. The third
author is a member of the Center for Mathematical Analysis,
Geometry and Dynamical Systems, IST  and was  partially supported
by FCT(Portugal) through program  POCTI/FEDER and Project
POCTI/MAT/44068/2002.}

\begin{document}

\begin{abstract}
We describe smooth rational projective algebraic surfaces $X$,
over an algebraically closed field of characteristic different
from $2$, having an even set of four disjoint $(-2)$-curves $N_1,\ldots,N_4$, i.e. 
such that $N_1+\cdots+N_4$ is divisible by 2 in $\Pic(X)$.
\end{abstract}

\maketitle

\section{Introduction}

Let $X$ be a smooth projective algebraic  surface over an
algebraically closed field  ${\mathbf k}$ of characteristic $\ne 2$. A
set of $\nu$ disjoint nodal (i.e $(-2)-$curves) curves $N_1,\ldots,N_{\nu}$ is called
an {\it even set} if there exists $L\in \Pic(X)$ such that
 $2L\equiv N_1+\cdots+N_{\nu}$. Since $K_XL=0$, $L^2$ is even and
therefore the number $\nu$ is divisible by
four because $4L^2=-2\nu$.

Assume that $X$ is rational, and that  $N_1,\ldots,N_{\nu}$ is an
even set of disjoint nodal curves  with  $\nu\geq 8$. Let
$\eta:X\to \Sigma$ be the map which contracts the curves $N_i$,
$i=1,\ldots,\nu$, to nodes.

By Theorem 3.2 of \cite{dmp}, there exists a fibration
$f:\Sigma\to \pp^1$ such that the general fibre of $F$ is a smooth
rational curve and having ${\nu}/{2}$ double fibres each
containing two nodes of $\Sigma$.

Equivalently there exists a fibration $g:X\to \pp^1$ with general
smooth rational fibres, having ${\nu}/{2}$ fibres, each containing
two nodal curves $N_k, N_l$ and of the form $N_k+2\Gamma+N_l$, where $\Gamma$ is a curve such that $\Gamma^2=-1$,
$K_X\Gamma=-1$. Since a rational fibration has no multiple fibres, by Zariski's lemma, (see \cite{bpv}),  such a curve  $\Gamma$ is necessarily 1-connected and as such  there is
a birational morphism contracting $\Gamma$ to a smooth point.
 Remark  that such a fibre can be obtained from a
minimal model of the fibration by first blowing up a point $x$,
then blowing up the double point of the total transform of the
ruling through $x$, and possibly further blowing-ups, cf.\ Example
1 in \cite{dmp}.

 Note also that any set of $4m$ of these nodal curves contained in $2m$ such
fibres is necessarily even.

By the proof of   Theorem  3.2 of \cite{dmp},  an even set
of four disjoint nodal curves  orthogonal to the fibres of
a fibration $g:X\to
\pp^1$, with smooth rational fibres, is necessarily  contained in two fibres of $g$ as above. 

Again  by Theorem 3.3 in \cite{dmp}, if a rational surface $Y$
with Picard number $\rho(Y)= 10-K_Y^2$ 
contains $\alpha$,  $\alpha\geq 3$,  disjoint nodal curves,
then $\alpha\leq \rho(Y)-2 $. Furthermore,  if  $\alpha=\rho(Y)-2 $,
then   $\alpha=2\beta $ is even, and $Y$ is obtained from a relatively minimal ruled
rational surface
$\F_e:= {\rm Proj}(\OO_{\pp^1}\oplus \OO_{\pp^1}(e))$, $e\ge 0$, with the process described above, i.e. by blowing up:
\begin{itemize} 
\item $\beta$ points $p_1, p_2, \ldots, p_{\beta}$ in distinct
fibres $F_1,F_2, \ldots,F_{\beta}$ of the same ruling of $\F_e$;
\item the point $q_i$, $i=1,2,\ldots,\beta$, which is the
intersection of the strict transform of $F_i$ with the exceptional
curve over $p_i$.
\end{itemize}
We keep the terminology introduced in \cite{dmp} and we will call
$Y$ as above the {\it standard example} of a rational surface with
$\rho(Y)-2$ disjoint nodal curves. For other terminology see ``Notation and conventions'' below.

In this note we complete the above results by proving the following:

\begin{thm}\label{TM} Let $X$ be a smooth rational surface
containing an even set of  4 disjoint nodal curves
$N_1,\ldots,N_4$. Then there exists a fibration $g:X\to \pp^1$
with smooth rational fibres, having $2$ fibres, each containing
two nodal curves $N_k, N_l$ and of the form $N_k+2\Gamma+N_l$,
where $\Gamma$ is a curve such that $\Gamma^2=-1$,
$K_X\Gamma=-1$.
\end{thm}

From  Theorem \ref{TM} we obtain:

\begin{cor} \label{1}
If $X$ is a smooth rational surface such that $X$ contains an even
set of four nodal curves, then there is a birational morphism $\pi:
X\to Y$, where $Y$ is a standard example with $K_Y^2=4$, which
maps the four nodal curves of $X$ to the four nodal curves of $Y$.
\end{cor}

\begin{proof}

Let $N_1,..., N_4$ be the nodal curves.
By Theorem \ref{TM}, $X$ has a fibration $ f:X\to \pp^1$ in
smooth rational curves having fibres $N_i+2\Gamma_1+N_j$, and
$N_k+2\Gamma_2+N_l$, $\{i,j,k,l\}=\{1,2,3,4\}$ and it is clear that, by contracting the
$(-1)$-curves in the other fibres and possibly some curves contained in $\Gamma_m$, $m=1,2$, we get to a standard example.
\end{proof}

\begin{cor}
Let $X$ be a weak Del Pezzo surface containing an even set of four
nodal curves.

Then there exists a birational morphism $\pi:X\to Y$, where
$K_Y^2=4$ and  $Y$ is obtained from $\F_a$, with $a=0$, $1$,
or $2$, by blowing up:
 \begin{itemize}

\item two points $p_1$, $p_2$ in distinct fibres $F_1,F_2$ of the
same ruling of $\F_a$;
\item the point $q_i$, $i=1,2$, which is the intersection of the
strict transform of $F_i$ with the exceptional curve of $p_i$.
\end{itemize}
Furthermore, in case $a=2$, none of the blown-up points lies on
the $(-2)$-curve of $\F_2$.
\end{cor}

\begin{proof}
Consider the birational morphism $\pi:X\to Y$ of
Corollary \ref{1}. Since $-K_X$ is nef, also $-K_Y$ has to be nef.
In particular, there are no rational curves $C$ on $Y$ or $X$ such that
$C^2\leq-3$. This implies the assertion.
\end{proof}

Our interest in these results arose in the course of investigating
surfaces of general type with an involution, cf.\ \cite{CCM}.
However it seems to us of independent interest.

The main ingredients used for the proof of Theorem \ref{TM}, which
is presented in section \ref{s:4}, are some facts on adjoint
systems on rational surfaces, which are collected in section
\ref{rat}, and  Lemma \ref{-1}, which is proved in section \ref{curve}.

\subsection*{Notation and conventions.} We work over any
algebraically closed field ${\mathbf k}$ of characteristic $\ne 2$.

All surfaces are projective algebraic varieties of dimension $2$
over ${\mathbf k}$. We do not distinguish between line bundles and divisors
on a smooth variety. Linear equivalence is denoted by $\equiv$ and
numerical equivalence over $\Q$ by $\sim$. The intersection
product of divisors (line bundles) $A$ and $B$ is denoted  by $A
B$.  As usual, given a divisor $D$ on a surface, $|D|$ will be the complete linear system of the effective divisors $D'\equiv D$.

By a  {\it curve} on a smooth surface $X$ we mean an effective, non zero
divisor on $X$. However a {\it $(-1)$-curve} (resp. {\it $(-2)$-curve} or {\it nodal curve}) is
an {\it irreducible} smooth rational curve $C$ such that $C^2=-1$ (resp.
$C^2=-2$).
 A {\it $(-1)$-divisor}  on $X$ will be a divisor $D$ satisfying $D^2=-1$
and $K_XD=-1$.

A smooth surface $X$ is called a {\it weak del Pezzo} surface if
$-K_X$ is big and nef. The remaining notation is standard in
algebraic geometry.

\section{Some properties  of rational surfaces}\label{rat}

In this section we list  some  properties of rational surfaces,
which we will  need later. The properties on adjoint systems
listed below can be also phrased in terms of Mori's theorem on the
cone (cf.\ \cite{reid}), but here, for the reader's convenience,
we state and prove them in the form we will need.

\begin{lem}\label{proper}
Let $X$ be a rational surface. Then:
\begin{enumerate} [(i)]

\item If $D$ is a $(-1)$-divisor, then either $|D|\neq \emptyset$ or $|K_X-D|\neq \emptyset$.

\item Assume $-K_X$ is nef and big. Then each effective
$(-1)$-divisor $D$ either contains a $(-1)$-curve or
$K_X^2=1$ and $D\in |-K_X+ A|=|-K_X|+A$, where $A$ is an effective divisor such that $K_XA=0, A^2=-2$.
\end{enumerate}
\end{lem}

\begin{proof}
The first assertion is an immediate consequence of the
Riemann-Roch theorem, because $\chi (\OX)=1$.

For (ii),  we note that, if $D$ is not irreducible, since $-K_X$
is nef and big, there is one component $\Delta$ of $D$ such that
$-K_X\Delta=1$ and every other component $\theta$ satisfies
$-K_X\theta=0$. By the index theorem  (see, e.g., Corollary 2.4 in
\cite{badescu}) and the adjunction formula, we see that $\theta^2=-2$ for each such $\theta$ and
either  $\Delta^2<0$ and $\Delta$ is a $(-1)$-curve or  $\Delta^2= K_X^2=1$ and $\Delta\sim -K_X$. In the last case,  $\Delta \equiv -K_X$ because numerical equivalence
 coincides with linear equivalence
 on rational surfaces. 
\end{proof}

\begin{lem}\label{adjunction1}
Let $D$ be a nef curve on a regular surface $X$ such that
$p_a(D)\geq 1$. If $K_X+D$ is not nef, then any irreducible curve
$\theta$ satisfying $\theta(K_X+D)<0 $ is a $(-1)$-curve $\theta$
such that $\theta D=0$.
\end{lem}

\begin{proof}
Since $X$ is regular, we have that $h^0(X, K_X+D)\geq p_a(D)\ge 1$.

Assume $K_X+D$ not nef. Then there is an irreducible curve $\theta$  such that $\theta(K_X+D)<0$. The curve $\theta$ is a component of the fixed part of $|K_X+D|$, and so   $\theta^2<0$.
 Since $D$ is nef, we have $K_X\theta<0$,
i.e.\ $\theta$ is a $(-1)$-curve and $\theta D=0$.
\end{proof}

\begin{lem}\label{adjunction2}
Let $D$ be a curve on a rational surface $X$ such that $p_a(D)= 1$ and
$D^2\geq 1$. If $K_X+D$ is nef, then $D\equiv-K_X$. In particular
$K_X^2\geq 1$.
\end{lem}

\begin{proof} As in the previous Lemma,  $h^0(X, K_X+D)\neq 0$. Since  $(K_X+D)D=0$, the index theorem and the hypothesis $K_X+D$ nef  imply that  $K_X+D\sim 0$, hence  $K_X+D\equiv 0$, because $X$ is rational.
\end{proof}

\begin{prop}\label{adjunction}
Let $D$ be a nef and big  2-connected curve  on a rational surface $X$ with
$K_X^2\le 0$. Assume that $p:=p_a(D)\geq 2$ and $K_X D<0$. If
$K_X+D$ is nef, then the following possibilities can occur:
\begin{enumerate} [(i)]

\item  $(K_X+D)^2=0$, and $K_X+D\equiv (p-1) G$, where $|G|$ is a
pencil of rational curves without base points such that $GD=2$;

\item  $(K_X+D)^2>0$, and the general curve $D_1$ in $|K_X+D|$ is
irreducible satisfying $ p_a(D_1)<p_a(D)$. 
\end{enumerate}

\end{prop}

\begin{proof}
Since  $X$ is regular, $h^0(X, K_X+D)=p\geq 2$.

Write $|K_X+D|=| M|+F$ where $|M|$ is the moving part and $F$ the
fixed part of the linear system $|K_X+D|$.
Since $D$ is 2-connected, $\omega_D$ has no base points (see, e.g.,
\cite{cfm}, Proposition A.7, or \cite{ma}) and so the nef divisor $D$ satisfies $DF=0$.  This implies
that, if $F \ne 0$, every curve $\theta$ contained in $F$ is such that $\theta D=0$
 and so, by the index theorem, $\theta^2<0$.
 In particular, if $F \neq 0$, then $MF>0$, because
$K_X+D$ is nef. Note that          $DM=D(M+F)=D(K_X+D)$ is even, hence the equality
$(K_X+D)M=M^2+MF$ implies, by the adjunction formula, that $MF$ is
even.

Suppose that the general curve $M$ in  $|M|$ is  reducible. Then $|M|$
is composed with a pencil $|G|$, $p \geq 3$,   $M\equiv (p-1)G$ and
$GD=2$. Note that $p \geq 3$ implies in particular that $D^2\geq
5$, because $K_XD<0$.  Since $GD=2$, the index theorem implies
that $G^2=0$. Note that  $GK_X=GF-GD=GF-2$ and thus $GF$ is even.
Now $((p-1)G) (K_X+D)=(p-1)GF$. Since
$$
((p-1)G) (K_X+D) \leq (K_X+D)^2= K_X^2+K_XD+2(p-1)<2(p-1),
$$
we conclude
that $GF=0$. So $F=0$ and we are in case (i).

Suppose now that the general curve $M$ in  $|M|$ is irreducible. We note that $h^0(X, M)=p$
and thus $h^0(M,\OO_M(M))=p-1$, because $X$ is regular.  Now note
that
$$K_XM=(K_X+D)M-MD=M^2+MF-2(p-1),$$
hence $K_XM+M^2=2M^2+MF-2(p-1)$. Since
$$
M^2+MF\leq (K_X+D)^2= K_X^2+K_XD+2(p-1)<2(p-1),
$$
one has $K_XM<0$ and therefore the series $ \OO_M(M)$ is non
special. By the Riemann-Roch theorem we obtain then
$$
p-1=M^2-M^2-{\frac {1}{2}  } MF+(p-1),
$$
hence $MF=0$ and thus $F=0$. If $M^2=0$, we have $p=2$ and we are
in case (i), whereas, if $M^2>0$, we are in case (ii).

In this case, since  $K_X^2\leq 0$ and $K_XD<0$, we have necessarily $p_a(M)<p_a(D)$.

\end{proof}

\section{Even sets of nodes and double covers} \label{curve}
Let $X$ be a smooth projective algebraic surface.
Given an even set of  disjoint nodal curves $N_1,\ldots,N_{\nu}$
on $X$, let $\pi\colon Y\to X$ be the double cover branched on
$N_1,\ldots,N_{\nu}$, defined by $2L\equiv N_1+\cdots+N_{\nu}$
(cf.\ pg.\ 42 in \cite{bpv}) and let $\eta\colon X\to \Si$, as
in the Introduction, be the map that contracts the curves $N_i$ to nodes. The
inverse image on $Y$ of a curve $N_i$ is a
$(-1)$-curve $\Delta_i$. Blowing these $(-1)$-curves down to
points $p_1,\ldots,p_{\nu}$, we obtain a smooth surface $\bar{Y}$
and a double cover $\bar{\pi}\colon \bar{Y}\to \Si$ branched
precisely over the singularities of $\Si$. Then we have  the
following commutative diagram:
\[
\begin{CD} Y@>  >>\bar{Y}\\ @V\pi VV@V\bar{\pi} VV\\ X@>\eta>>\Sigma,
\end{CD}
\]
Note that $\Sigma$ has canonical singularities, so that $K_\Sigma$
is a Cartier divisor. Moreover
$\bar\pi^*(K_\Sigma)=K_{\bar{Y}}$. Hence
\begin{align*}
& \kappa(Y) = \kappa(\bar{Y})=\kappa(X), & & K^2_Y=2K^2_X-{\nu}, \\
& K_{\bar{Y}}^2 = 2K^2_X,
 & & \chi(Y,\OO_Y) = \chi(\bar{Y},\OO_{\bar{Y}}) = 2-\frac{\nu}{4}.
\end{align*}

Finally we will need the following:

\begin{lem}\label{-1}  Let $\bar{E}$ be a
$(-1)$-curve of $\bar{Y}$ and let $E$ be the strict transform of
$\bar{E}$ in $Y$. Then   $ E$ is a component of ${\pi}^*(C)$ where
$C$ is an irreducible curve such that $K_XC=-1$, and such that,
for each nodal curve $N_i$, either $CN_i=2$ or $CN_i=0$.
\end{lem}

\begin{proof}
Since  $\bar{E}$ is a smooth curve,  $E$ meets each  of the
$(-1)$-curves $\Delta_i$ transversally in at most one point. Note
that $\sum \Delta_i\equiv \pi^*(L)$. Let $m$ be the number of curves
$\Delta_i$ having non-empty intersection with $E$. Then $E^2=-1-m$
and $K_YE=m-1$. Since $K_Y\equiv \pi^*(K_X+L)$,  we conclude that
$\pi^*(K_X)E=-1$.  Since the map $H^2(X,\Q)\to H^2 (Y,\Q)$ induced
by $\pi$ multiplies the intersection form by $2$, we conclude that
the curve $E$ is not invariant under the involution $\iota$ of $Y$
associated to $\pi$. Then, if $C=\pi(E)$, $\pi^*(C)=E+\iota(E)$
and $C$ is as stated.
\end{proof}

\section{The proof of Theorem \ref{TM}} \label{s:4}

We use the notation of the statement of Theorem \ref{TM}  and we
denote again by $L$ the line bundle such that
$N_1+\cdots+N_4\equiv 2L$. The line bundle $L$ satisfies $L^2=-2$,
 $K_XL=0$ and $|L|=\emptyset$.

We will need the following:

\begin{lem}\label{fibre}
If there exists a $(-1)$-curve $E$ such that $EL=1$
and $E$ meets transversally exactly two of the nodal curves, say
$N_1$, $N_2$, then $X$ is as in Theorem \ref{TM}.
\end{lem}

\begin{proof}
Since $EL=1$, $E+L$ is a $(-1)$-divisor and therefore by Lemma
\ref{proper}, (i), either $|E+L|\neq \emptyset$ or $|K_X-(E+L)|\neq \emptyset$.           The second possibility clearly does not occur, since,
otherwise, $2K_X$ would be effective. Therefore $|E+L|\neq \emptyset$.
 Since $N_3L=N_4L=-1$, we can write $E+L\equiv \Gamma+N_3+N_4$
where $\Gamma$ is an effective $(-1)$-divisor.

Note that $E(E+L)=0$ implies $E\Gamma=0$ and actually
$E\cap\Gamma=\emptyset$.  In fact otherwise $E$ would be a
component of $\Gamma$, hence $E+L\equiv E+\Delta$, where $\Delta$ is an
effective divisor implying that  $|L|\neq \emptyset$.

By the Riemann-Roch theorem, $h^0(X, 2E+N_1+N_2)\geq 2$. Now, the
relation  $2E+N_1+N_2+N_3+N_4\equiv 2(E+L)\equiv2\Gamma+2N_3+2N_4$
implies $2E+N_1+N_2\equiv 2\Gamma+N_3+N_4$. So $|2E+N_1+N_2|$ is a
pencil of rational curves without base points having fibres as in
the statement.
\end{proof}

Now we can give the:

\begin{proof}[Proof of Theorem \ref{TM}]
Since, by contracting $(-1)$-curves  disjoint from $N_1,$
$ \ldots, N_4$,
  we still obtain a surface having an even
set of 4 disjoint nodal curves, we will from now on make the
following:

\begin{assu} \label{EL>0}
There is no  $(-1)$-curve on $X$ disjoint from the curves
 $N_1, \ldots, N_4$,
i.e., for every  $(-1)$-curve $E$, one has $EL\ge1$.
\end{assu}

We will argue by contradiction. So suppose that there is no
fibration as in the statement. This implies that $K_X^2< 4$, by
Theorem 3.3 in \cite{dmp}. Moreover, by Lemma \ref{fibre} and
Assumption \ref{EL>0}, only the following two cases are possible:
\begin{enumerate}[ (I)]
\item given a $(-1)$-curve $E$, one has $EL\geq 2$;
thus by the index theorem $K^2_X\le1$, and $K_X^2=1$ if and only
if $EL=2$ and $-K_X\equiv E+L$; or \label{EL>1}
\item given a $(-1)$-curve $E$, one has $EL=1$ and
$E$ intersects exactly one of the nodal curves, say $EN_1=2$; thus
by the index theorem $K^2_X\le1$ and $K^2_X=1$ if and only if
$-K_X\equiv E+N_1$. \label{EL=1}
\end{enumerate}

  The surface $\bar{Y}$ as in Section \ref{curve} is not minimal, because $K_{\bar{Y}}^2<8$ and so  
there exists on $X$  an irreducible curve $C$ as in
Lemma \ref{-1} meeting $m$ of the nodal curves, say
$N_1,\ldots,N_m$ with $m\geq 0$,  and satisfying $CL=m$. Since $C$
is irreducible and $K_XC=-1$, one has $C^2\geq -1$. Furthermore
Assumption \ref{EL>0} means that, if $C^2=-1$, then necessarily
$m>0$.

Set $D:=C+N_1+\cdots+N_m$. Remark that
$$DN_i=0, i=1,\ldots, 4, \quad\text{and therefore } \quad DL=0 \text{
and } D^2=C^2+2m.$$
The curve $D$ is nef, big and 2-connected. Since $K_XD=-1$, the index theorem yields
$K_X^2\leq1$, and
$$K^2_X=1 \Longleftrightarrow -K_X\equiv D.$$

We start by considering this case.

\subsection{The case $K_X^2=1$.}
Since $-K_X\equiv D$, $D^2=1$ and so either $D=C$ is irreducible or  $C$ is a $(-1)$-curve  intersecting only one of the nodal curves, say $N_1$, and $D=C+N_1$. We notice that in both
cases $|-K_X|= |D|$ is a pencil without fixed components.

Since $K_X^2=1$,  $-K_X+L$ is a $(-1)$-divisor, and,  by
Lemma \ref{proper}, (i), $|-K_X+L|\neq \emptyset$   because $-K_X(2K_X-L)<0$. Since, for each
$i=1,\ldots,4$, $N_i(-K_X+L)=-1$, we can write
$-K_X+L\equiv \Gamma+N_1+N_2+N_3+N_4$ where $\Gamma$ is  an
effective $(-1)$-divisor.  Note that
\begin{equation} \label{Gamma+L=-K}
\Gamma+L\equiv-K_X.
\end{equation}

\subsubsection{Claim: $\Gamma$ is irreducible.}
\label{claim:Gamma_irr} Since $-K_X$ is nef and big, by  Lemma
\ref{proper}, (ii), $\Gamma$ contains one $(-1)$-curve
$\gamma$ which will satisfy one of the cases \eqref{EL>1} or  \eqref{EL=1}.

If $\gamma L\ge2$, then, by  case  \eqref{EL>1} and identity \eqref{Gamma+L=-K}, we
have $-K_X\equiv \gamma+L\equiv\Gamma+L$, which implies
$\gamma=\Gamma$, i.e.\ the claim.

If $\gamma L=1$, then, again by \eqref{Gamma+L=-K}, one has
$\gamma\Gamma=0$. Since $\gamma$ is in case  \eqref{EL=1}, there is a nodal curve $N_i$
such that $\gamma N_i=2$ and $\gamma+N_i\equiv -K_X$. Since
$N_i\Gamma=N_i(-K_X-L)=1$ by \eqref{Gamma+L=-K}, then
$N_i\gamma=2$ implies that $N_i$ is also a component of $\Gamma$.
But, always by identity  \eqref{Gamma+L=-K}, one then has $-K_X\equiv
\gamma+N_i\leq\Gamma\equiv -K_X-L$, which implies $-L\geq 0 $, that
is impossible.

\subsubsection{Claim: every  $(-1)$-curve $E\neq\Gamma$
satisfies $EL=1$, hence there is a nodal curve $N_i$ such that
$EN_i=2$ and $-K_X\equiv E+N_i$.} \label{claim:EL=1}

Suppose that $E$ is a $(-1)$-curve such that $EL=2$.
Then, by \eqref{EL>1} and \eqref{Gamma+L=-K}, one has $E+L\equiv
-K_X$, hence $E=\Gamma$. The last assertions follow by
case \eqref{EL=1}.

\subsubsection{Claim: there are  $(-1)$-curves $E_1,
E_2$ different from $\Gamma$, such that $E_1E_2=1$,
$E_1N_1=E_2N_2=2$ and $-K_X\equiv E_1+N_1\equiv E_2+N_2$.}

Since $-K_X$ moves in a pencil without fixed components and
$-K_XN_1=0$,  there is a curve $E_1+N_1$ in the pencil $|-K_X|$,
where $E_1$ is an effective $(-1)$-divisor. 

The curve $E_1$ is irreducible.  Indeed, by Lemma \ref{proper},
(ii), $E_1$ contains a $(-1)$-curve $\theta$. Remark
that $\theta\neq\Gamma$, otherwise by \eqref{Gamma+L=-K} we would
have $\theta+L\equiv E_1+N_1$, which would imply $L>0$, a
contradiction. Hence, by Claim \ref{claim:EL=1}, there exists one
of the nodal curves $N_i$ such that $-K_X\equiv \theta+N_i\equiv
E_1+N_1$. Since $\theta\leq E_1$ and $|-K_X|$ has no fixed
components, this implies $\theta=E_1$.

The curve $E_2$ is found by applying the same reasoning to the
fibre of the pencil $|-K_X|$ which contains $N_2$.

Since $-K_X\equiv E_1+N_1\equiv E_2+N_2$, for both curves $E_1$,
$E_2$ we are in case \eqref{EL=1}, and not case \eqref{EL>1}. Therefore
$E_1N_2=E_2N_1=0$, which implies that $E_1E_2=1$.

\subsubsection{Claim: the linear system $|E_1+E_2|$ is a base point
free pencil of rational curves, the curve $2\Gamma+N_3+N_4$ sits
in the pencil $|E_1+E_2|$, which has at least three reducible
fibres.}

Notice that $-2K_X\equiv E_1+N_1+E_2+N_2\equiv
2\Gamma+N_1+N_2+N_3+N_4$, whence the first two assertions follow.
For the last assertion, remark that  $\rho(X)=9$, thus
$|E_1+E_2|$ contains yet another reducible fibre.

\bigskip
Now we can conclude the proof for the case $K^2_X=1$.

A reducible fibre of $|E_1+E_2|$  contains at least one 
$(-1)$-curve $G$.  So there is a $(-1)$-curve $G$ such that
$GE_1=GE_2=G\Gamma=0$. Since $G\neq \Gamma$, one has $GL=1$ by
Claim \ref{claim:Gamma_irr} and so $G$ is in  case  \eqref{EL=1}.  On the other hand, $1=-K_X
G=G(E_i+N_i)=GN_i$, $i=1,2$, which is not possible in case \eqref{EL=1}.

\subsection{The case $K_X^2<1$.}

We start with the following:

\subsubsection{Claim: every $(-1)$-curve $E$ satisfies
$EL\geq 2$, i.e.\  we are in case \eqref{EL>1}.} \label{claim:caseEL>1}

Suppose otherwise, namely suppose there is a
$(-1)$-curve $E$  for which case \eqref{EL=1} holds, i.e.\ $EL=1$,
$EN_1=2$ and $EN_i=0$, $i=2,3,4$. Hence the curve $A:=E+N_1$ is
nef,  $p_a(A)=1$ and   $AL=0$. Since $(K_X+A)^2<0$, then
$K_X+A$ is not nef and so, by Lemma \ref{adjunction1}, there
exists a $(-1)$-curve $\theta$ such that $\theta
A=0$. Then one has $(\theta+L) A=0$ and, therefore,
$(\theta+L)^2<0$ by the index theorem. This implies $\theta L=1$,
namely $\theta$ is as in case \eqref{EL=1}, i.e.\ there is a nodal
curve, say $N_2$, such that $N_2\theta=2$. But then
$(N_2+\theta)^2=1$ and $A(N_2+\theta)=0$, which contradicts the
index theorem. This proves the claim.

\medskip

Now we consider again the nef and big 2-connected  curve $D:=C+N_1+...+N_m$, which satisfies $K_X D=-1$. In particular $p_a(D)\geq 1$.

\subsubsection{Claim: $K_X+D$ is nef and moreover $D^2\geq 3$, $p_a(D)\geq 2$} \label{claim:K+Dnef}

Suppose  that $K_X+D$ is not nef. By Lemma \ref{adjunction1} there is a $(-1)$-curve $E$ such that $DE=0$. By Claim
\ref{claim:caseEL>1}, one has $(E+L)^2>0$. Since $DL=0$, one
also has $D(E+L)=0$. This gives a contradiction to the index
theorem and so $K_X+D$ is nef.

 In particular $0\leq (K_X+D)^2=K_X^2+2K_XD+D^2$. Since $K_X^2\leq 0$ and  $K_XD=-1$, we obtain $D^2\geq 2$. Since $D^2$ is odd by the adjunction formula, we have proved the last two assertions.

\subsubsection{Claim: there is a positive dimensional linear system $|M|$ whose general
curve $M$ is irreducible, smooth, rational and such that $ML=0$.
} \label{claim:last}
We note first that $mK_X+D$ is orthogonal to $L$, for any $m\in \N$.

If $(K_X+D)^2=0$,
by Lemma \ref{adjunction} one has   $K_X+D\equiv (p-1) G$, where
 $|G|$ is a pencil
 of rational curves without base points and we have proven the claim.
If $(K_X+D)^2>0$, 
 again by Lemma \ref{adjunction} 
 the general curve $D_1
\in |K_X+D|$ is irreducible with $ p_a(D_1)<p_a(D)$. If $p_a(D_1)=0$,
 again we proved the claim.

 If $p_a(D_1)>0$ notice that, since $D_1$ is orthogonal to
$L$ and $D_1^2>0$, we   can show as in Claim  \ref{claim:K+Dnef}  that $K_X+D_1\equiv 2K_X+D$ is nef. So,   by Lemma \ref{adjunction2}, $p_a(D_1)\geq 2$  and, as in the previous paragraph, by Lemma \ref{adjunction},   either we find a linear system
as in the claim or the general curve $D_2 \in |K_X+D_1|$ is irreducible
and satisfies $0< p_a(D_2)<p_a(D_1)<p_a(D)$.

It is clear that by iterating this procedure we eventually find a linear system as in the claim.

\medskip
Now  we can finish our proof of this case, and therefore of the
theorem.

 Consider   the    positive dimensional linear system
$|M|$, whose existence is proved in Claim \ref{claim:last}, and let $M$ be a general curve of $|M|$.

If $M^2=0$, then $|M|$ is a base-point-free pencil. Then by
\cite{dmp}, as recalled in the introduction, the pencil $|M|$ is
as in the statement of Theorem \ref{TM}, contradicting our
assumption.

If $M^2>0$, we also have a contradiction. Indeed, since  $M$ is
rational and smooth, $L$ is trivial on $M$ and thus, if $\pi:Y\to
X$ is the double cover branched on $N_1+\cdots+N_4$, one has
$\pi^*(M)=M_1+M_2$, where $M_1M_2=0$ and $M_i^2=M^2$. Since
$M^2>0$, this contradicts the index theorem.
\end{proof}

\bigskip

\bigskip

\bigskip

\noindent {\it Authors' adresses}:

\medskip

\small{ \noindent Alberto Calabri, Dipartimento di Matematica,  Universit\`a degli Studi di Bologna, Piazza di Porta San Donato, 5,  I-40126 Bologna, Italy.

\noindent {\it e-mail}: calabri@dm.unibo.it

\medskip

\noindent Ciro Ciliberto, Dipartimento di Matematica,  Universit\`a
degli Studi di Roma ``Tor Vergata'', Via della Ricerca Scientifica, I-00133 Roma, Italy

\noindent {\it e-mail}:  cilibert@mat.uniroma2.it

\medskip

\noindent Margarida Mendes Lopes,
Departamento de  Matem\'atica, Instituto Superior T\'ecnico,  
Universidade T\'ecnica de Lisboa, 
Av.~Rovisco Pais, 1049-001 Lisboa, Portugal

\noindent {\it e-mail}: mmlopes@math.ist.utl.pt
 }

\end{document}